\documentclass[11pt,twoside]{article}
\usepackage{latexsym, amssymb, amsthm, amsmath, epsfig, setspace, tikz, slashed,mathrsfs}
\usepackage{mathrsfs}
\usepackage[toc,page]{appendix}
\usepackage[all]{xy}
\usepackage[retainorgcmds]{IEEEtrantools}
\usetikzlibrary{decorations.markings}
 \usetikzlibrary{arrows}
\usepackage{hyperref}
\usepackage[margin=10pt,font=small,labelfont=bf]{caption}

\theoremstyle{plain}
\newtheorem{theorem}{Theorem}[section]
\newtheorem{proposition}[theorem]{Proposition}

\theoremstyle{definition}

\newtheorem{corollary}[theorem]{Corollary}

\newtheorem{ex}[theorem]{Example}

\newenvironment{example}%
{\begin{ex}}%
{\hfill$\blacksquare$\end{ex}}

\newenvironment{renumerate}%
{%
\begin{enumerate}}%
{\end{enumerate}%
}%

\newenvironment{remark}%
{\vskip6pt%
\noindent%
{\it Remark.}}%
{\vskip6pt}

{\vskip6pt%
\noindent%
{\it Remarks}. %
\begin{renumerate}}%
{\end{renumerate}\vskip6pt}

\makeatletter
\def\Ddots{\mathinner{\mkern1mu\raise\p@
\vbox{\kern7\p@\hbox{.}}\mkern2mu
\raise4\p@\hbox{.}\mkern2mu\raise7\p@\hbox{.}\mkern1mu}}
\makeatother

\newcommand{\C}{\text{$\mathbb C$}}

\newcommand{\Z}{\text{$\mathbb Z$}}

\newcommand{\T}{\text{$\mathbb{T}$}}
\newcommand{\TM}{\text{$\mathbb{T}M$}}

\renewcommand{\tilde}{\widetilde}

\newcommand{\J}{\text{$\mathcal{J}$}}

\newcommand{\M}{\text{$\mathcal{M}$}}

\newcommand{\G}{\text{$\mathcal{G}$}}

\newcommand{\Ch}{\mathrm{Ch}}

\newcommand{\gf}{\text{$\varphi$}}

\newcommand{\tensor}{\otimes}

\newcommand{\mc}[1]{\text{$\mathcal{#1}$}}
\newcommand{\msc}[1]{\text{$\mathscr{#1}$}}

\newcommand{\into}{\longrightarrow}
\newcommand{\noqed}{\let\qed\relax}

\newcommand{\Diff}{\mathrm{Diff}}

\newcommand{\gcs}{generalized complex structure}

\newcommand{\gcss}{generalized complex structures}

\newcommand{\gc}{generalized complex}

\date{} \usepackage{color} \definecolor{tocolor}{rgb}{.1,.1,.5}
\definecolor{urlcolor}{rgb}{.2,.2,.6}
\definecolor{linkcolor}{rgb}{.1,.1,.6}
\definecolor{citecolor}{rgb}{.6,.2,.1}
\hypersetup{backref=true, colorlinks=true, urlcolor=urlcolor,
  linkcolor=linkcolor, citecolor=citecolor}

\begingroup
\hypersetup{linkcolor=blue}
\endgroup

\numberwithin{equation}{section}

\makeatletter
\newcommand{\subjclass}[2][2010]{%
  \let\@oldtitle\@title%
  \gdef\@title{\@oldtitle\footnotetext{#1 \emph{Mathematics subject classification.} #2}}%
}
\newcommand{\keywords}[1]{%
  \let\@@oldtitle\@title%
  \gdef\@title{\@@oldtitle\footnotetext{\emph{Keywords.} #1.}}%
}
\makeatother

\begin{document}
\title{A remark on the number of components of the space of \gcss}
\author{Gil R. Cavalcanti\thanks{{\tt gil.cavalcanti@gmail.com}} \\
       Department of Mathematics\\
Utrecht University\\
}

\maketitle

\abstract{We give examples of generalized complex four-manifolds whose moduli space has infinitely many components.}



\section{Complex structures and Chern classes}

Given a manifold $M$, there are three groups of symmetries one can consider, from biggest to smallest: $\mathrm{Diff}(M)$, the group of all diffeomorphisms of $M$;  $\mathrm{Diff}_0(M)$, the group of diffeomorphisms which induce the identity map in cohomology and $\tilde{\mathrm{Diff}_0}(M)$, the connected component of the identity.  Accordingly, given a geometric structure, $\mc{S}$, one has three moduli ``spaces" on $M$: let $\msc{S}$ be the space of all $\mc{S}$-structures on $M$ and define the sets, from smallest to biggest:
$$\mc{M}_{\mc{S}} = \frac{\msc{S}}{\mathrm{Diff}(M)},\qquad \mc{M}^0_{\mc{S}}= \frac{\msc{S}}{\mathrm{Diff}_0(M)},\qquad \tilde{\mc{M}^0_{\mc{S}}} = \frac{\msc{S}}{\tilde{\mathrm{Diff}_0}(M)}.$$
One of the most basic things one can ask about these spaces is their cardinality. In particular, if the spaces are already known to be enumerable, one would like to know at least if they are finite or not.

We will be interested in generalized complex structures, but lets first comment on the case of complex structures. We let $\msc{I}$ be the space of all complex structures and $\msc{I}_{alm}$ be the space of all almost complex structures on $M$. Kuranishi's theorem \cite{MR0176496} implies that $\M^0_{\msc{I}}$,  is locally connected and hence it can have at most enumerably many components. There is a natural inclusion $\msc{I} \hookrightarrow \msc{I}_{alm}$ and we will study the question of finiteness  of components of the moduli ``spaces" above via the Chern polynomial:
$$\Ch:\msc{I}_{alm} \into H^{\bullet}(M;\Z).$$
The point being that if the restriction of the Chern polynomial to $\msc{I}$ has infinite image, then $\mc{M}^0_{\msc{I}}$ has infinitely many components, while if $\Ch(\msc{I}_{alm})$ is finite, this simple topological invariant can not be used to determine if $\mc{M}^0_{\msc{I}}$  is infinite or not. Of course, one can use the same tool to tackle the smaller space, $\mc{M}_{\msc{I}}$, but now the Chern polynomial takes values in the quotient space $H^{\bullet}(M;\Z)/\mathrm{Diff}(M)$.

The simplest situation in which one can expect to use the Chern polynomial effectively is that of an almost complex four manifold, $M^4$. In this case, there are only two Chern classes and we have the relations
\begin{equation}\label{eq:relations}
\begin{aligned}
c_2(M) & = \chi(M);\\
c_1^2(M) &= 3\sigma(M) + 2\chi(M);\\
c_1 &= w_2 \mbox{ mod } 2,
\end{aligned}
\end{equation}
where $\chi$  is the Euler class,   $\sigma$, the signature and $w_2$ the second Stiefel--Whitney class. These show that the topology of $M$ determines $c_2$ and constrains $c_1$. In particular the question of whether of not the Chern polynomial has finite image boils down to which values of $c_1$ one can achieve.

\begin{example}[Almost complex structures on $\C P^2$]\label{ex:cp2v1}
Since $H^2(\C P^2;\Z) = \Z$, with generator $a$ satisfying $a^2([\C P^2])= 1$, we have that $c_1(\C P^2) = \gamma a$ for some integer $\gamma$. Therefore, using \eqref{eq:relations} we have $\gamma^2 = 9$ and hence $ c_1 = \pm  3 a$ and the Chern polynomial of $\C P^2$ is
\begin{equation}\label{eq:cp2 chern polynomial}
\Ch(\C P ^2) = 1 \pm 3 a + 3 a^2;
\end{equation}
that is, the image of  $\Ch:\tilde{\mc{M}^0_{\msc{I}_{alm}}}\into H^\bullet(M;Z)$ only has two points.

Since complex conjugation of $\C P^2$ maps $a$ to $-a$ we have that, for $\C P^2$,
$$\Ch:\mc{M}_{\msc{I}_{alm}}\into  H^\bullet(M;Z)/\Diff(M)$$
is the constant map and this extreme example shows that due to purely topological reasons Chern classes do not help to determine the number of components of $\mc{M}_{\msc{I}}$.
\end{example}

This example is a particular case of  a more general finiteness result:

\begin{proposition}\label{prop:definite}
Let $M$ be a four-manifold with definite intersection form. Then the Chern polynomial
$$\Ch: \M^0_{\msc{I}_{alm}} \into H^{ev}(M;\Z);\qquad [\mc{I}]{\longmapsto} \Ch(\TM),$$
has finite image. 
\end{proposition}
\begin{proof}
Indeed, if $\{a_1,\cdots a_n\}$ is an integral basis of $H^2(M;\Z)/\mathrm{Tor}(H^2(M;\Z))$, then, modulo torsion, $c_1 = \sum \gamma_i a_i$ and 
$$\|3 \sigma(M) + 2 \chi(M)\| =  \|c_1^2([M])\| = \|\sum \gamma_i^2 a_i^2(M)\| \geq \sum \gamma_i^2,$$
showing that the vector $(\gamma_1,\cdots, \gamma_m) \in \Z^m$ lies in the ball or radius $\sqrt{\|3 \sigma(M) + 2 \chi(M)\|}$.  Since $\mathrm{Tor}(H^2(M;\Z))$ is also finite, we conclude that there are only finitely many polynomials which are the Chern polynomial of an almost complex structure.
\end{proof}

\section{Generalized complex structures and Chern classes}

The discussion above can be transported to the realm of \gc\ geometry. Indeed, Gualtieri proved a deformation theorem {\it \`a la} Kuranishi \cite{MR2811595} which automatically implies that the moduli space $\mc{M}_{\msc{J}}^0$ has at most countably many components. Since a \gcs\ is a complex structure on $\TM = TM \oplus T^*M$, we also have a corresponding set of Chern classes and can form the Chern polynomial:
$$\Ch:\M_{\msc{J}}^0 \into H^{\bullet}(M;\Z); \qquad [\J] {\longmapsto} \Ch(\TM).$$
If $\J$ is induced by a complex structure, then its $+i$-eigenbundle is $\T^{1,0}M = T^{0,1}M\oplus T^{*1,0}M$ and both summands are isomorphic, as complex vector bundles, hence the Chern polynomial of $\TM$ is the square of the Chern polynomial of $T^{*1,0}M$. If $\J$ is symplectic, then $\T^{1,0}M$ is isomorphic to  $T_\C M$, the complexification of $TM$, and the Chern polynomial of these bundles agree.

Besides the Chern classes of $\TM$, a \gcs\ $\J$ determines and is determined by a line subbundle $K_{\J} \subset \wedge T^*_\C M$, its {\it canonical bundle} and hence we have a corresponding first Chern class. For complex structures the canonical bundle agrees with the homonymous bundle from complex geometry and hence the first Chern class of $K_\J$ is the first Chen class of $T^{*1,0}M$. 
For symplectic structures $K_\J$ is the line generated by the form $e^{i\omega}$ and hence $c_1(K_{\J})=0$, as $K_\J$ has a nowhere vanishing section. Another class of interest is that of \gcss\ with nondegenerate type-change locus. In this case, the structure is generically symplectic but has more exotic behaviour along an anticanonical divisor (an embedded submanifold, $\Sigma$, of codimension two) and we have
\begin{equation}\label{eq:c1}
PD([\Sigma]) = - c_1(K_\J).
\end{equation}

Using Clifford action we get an isomorphism $\wedge^{top}\T^{1,0}M \tensor \overline{K} \cong K$ and hence the Chern classes are related by
$$c_1(\TM) = 2 c_1(K).$$

Further, given a \gcs\ $\J_1$, we can always pick a compatible metric $\G$ to obtain a second almost \gcs\ $\J_2 = \G \J_1$ and almost complex structures $I_\pm$. The Chern classes of $\TM$, $K_{\J_1}$, $K_{\J_2}$ and  $T^{1,0}_\pm M$,  the $+i$-eigenbundle of $I_\pm$. The relations between the bundles
$$\T^{1,0}_{\J_1}M \cong T^{1,0}_+M \oplus T^{1,0}_-M\qquad \mbox{and} \qquad \T^{1,0}_{\J_2}M  \cong T^{1,0}_+M \oplus T^{0,1}_-M,$$
give rise to relations between the Chern classes:
\begin{equation}\label{eq:chern polynomials}
c(\T^{1,0}_{\J_1}M )= c(T^{1,0}_+M)\cup c(T^{1,0}_-M)\qquad \mbox{and} \qquad c(\T^{1,0}_{\J_2}M) = c(T^{1,0}_+M) \cup c( T^{0,1}_-M).
\end{equation}

\begin{example}[Almost generalized complex structures on $\C P^2$]\label{ex:cp2v2}
There are two well known \gcss\ on $\C P^2$: the symplectic, $\J_\omega$, and the complex, $\J_I$. For the former
$$\Ch(\TM) = \Ch(T_\C M) = \Ch(T^{1,0}M)\cup \Ch(T^{0,1}M)  = (1+a)^3\cup (1-a)^3 = 1 -3 a^2$$
while for the latter
$$\Ch(\TM) = \Ch(T^{*1,0})^2 = (1+a)^6 = 1 - 6 a + 15 a^2, $$
where $a\in H^2(\C P^2;\Z)$ is a generator.

Modulo the action of $\Diff(\C P^2)$, these are the only possible Chern polynomials for any almost \gcs\ on $\C P^2$. Indeed, due to Example \ref{ex:cp2v1}  we know that the Chern polynomial of the almost complex structures $I_\pm$ must be of the form \eqref{eq:cp2 chern polynomial} and then the relations \eqref{eq:chern polynomials} give two possibilities for $c(\TM)$ according to whether the signs of the first Chern classes of $I_+$ and $I_-$ agree or not. If they have opposite signs, the Chern polynomial agrees with that corresponding to the symplectic structure, if they have the same sign, it agrees with that corresponding to the symplectic structure (possibly after complex conjugation of  $\C P^2$). So the Chern polynomial allows us determine the existence of  two components of the moduli space of \gcss\ on $\C P^2$, but no further.
\end{example}

\begin{remark}
Recently Goto and Hayano \cite{goto-hayano} and Torres and Yazinski \cite{torres-yazinski} produced examples of \gcss\ on $\C P^2$ whose complex locus has an arbitary number of components. The example above shows that Chern classes alone do no help to determine if those are in the same component of the moduli space or not. 
\end{remark}

As before we have a finiteness result for manifolds with definite intersection form:

\begin{corollary}
Let $M$ be a four manifold with positive definite intersection form. Then the Chern polynomial
$$\Ch: \M_{\msc{J}_{alm}}^0 \into H^{ev}(M;\Z);\qquad [\J]{\longmapsto} \Ch(\TM),$$
has finitely many points in its image. 
\end{corollary}
\begin{proof}
Due to Proposition \ref{prop:definite}, $\Ch(\mc{M}_{\msc{I}_{alm}}^0)$ is finite and due to \eqref{eq:chern polynomials} points in  $\Ch(\M_{\msc{J}_{alm}}^0)$ are determined by pairs of points in $\Ch(\mc{M}_{\msc{I}_{alm}}^0)$.
\end{proof}

In contrast, if the intersection form on $M^4$ is nondegenerate, the image of the Chern polynomial may contain infinitely many points and hence $\M_{\msc{J}}$ has infinitely many components.

\begin{theorem}[Elliptic surfaces]
Let $M$ be an elliptic surface with positive Euler characteristic. Then
$$\Ch: \M_{\msc{J}} \into H^{ev}(M;\Z)$$
has infinite image.
\end{theorem}
\begin{proof}
Under the hypothesis, $M$ has a symplectic structure, $\omega$ for which the fibers are symplectic and the fibration has at least one fishtail fiber and hence, after deformation, at least twelve. Then Goto and Hayano \cite{goto-hayano} (see also Cavalcanti and Gualtieri \cite{MR2312048}) proved that one can perform a multiplicity one logarithmic transform on a regular fiber to create a new manifold
$$\tilde{M} = M\backslash \pi^{-1}(D^2) \cup_\gf D^2 \times T^2,$$
and that $\tilde{M}$ admits a \gcs, $\J$, for which $\Sigma = \{0\}\times T^2 \subset D^2 \times T^2$ is a nondegenerate type-change locus. Since $M$ has enough fishtail fibers, $\tilde{M}$ is in fact diffeomorphic to $M$ and the type change locus, $\Sigma$,  is cohomologous to the fibers of the fibration. Hence we conclude that $c_1(K_\J) = - PD[\Sigma] \neq 0$. Since for the symplectic structure $c_1(K_{\J_\omega})=0$, we conclude that these structures are in different components of $\M_{\msc{J}}^0$.

One does not have to stop at one regular fiber and performing multiplicity-one logarithmic transforms at $k$ regular fibers, one produces a \gcs\ $\J$ with nondegenerate type change locus which represents $k [F]$ where $F$ is any regular fiber. Since, for different values of $k$,  the classes $kPD([F])$ are not in the same orbit of $\Diff(M)$ we conclude that $\Ch(\M_\msc{J})$ is infinite.
\end{proof}

\bibliographystyle{hyperamsplain}
\bibliography{references}

\end{document}